\documentclass{gtmon_a}
\pdfoutput=1
\usepackage{pinlabel}
\usepackage[all]{xy}
\input{diagxy}

\proceedingstitle{Exotic homology manifolds (Oberwolfach 2003)}
\conferencestart{29 June 2003}
\conferenceend{5 July 2003}
\conferencename{Workshop on Exotic Homology Manifolds}
\conferencelocation{Oberwolfach Mathematics Institute, Oberwolfach,
Germany}

\editor{Frank Quinn}
\givenname{Frank}
\surname{Quinn}

\editor{Andrew Ranicki}
\givenname{Andrew}
\surname{Ranicki}

\title{Stability in controlled $L$--theory}

\author{Erik Kj{\ae}r Pedersen}
\givenname{Erik Kj{\ae}r}
\surname{Pedersen}
\address{Department of Mathematical Sciences\\SUNY at Binghamton\\\newline
Binghamton\\New York 13901\\USA}
\email{erik@math.binghamton.edu}
\urladdr{}

\author{Masayuki Yamasaki}
\givenname{Masayuki}
\surname{Yamasaki}
\address{Department of Applied Science\\Okayama University of Science\\\newline
Okayama\\Okayama 700-0005\\Japan}
\email{masayuki@mdas.ous.ac.jp}
\urladdr{}

\volumenumber{9}
\issuenumber{}
\publicationyear{2006}
\papernumber{5}
\startpage{67}
\endpage{86}

\doi{}
\MR{}
\Zbl{}

\arxivreference{math.GT/0402218}

\keyword{controlled $L$--groups}
\subject{primary}{msc2000}{18F25}
\subject{secondary}{msc2000}{57R67}

\received{13 February 2004}
\revised{}
\accepted{13 February 2004}
\published{22 April 2006}
\publishedonline{22 April 2006}
\proposed{}
\seconded{}
\corresponding{}
\version{}

\makeatletter
\def\cnewtheorem#1[#2]#3{\newtheorem{#1}{#3}
\expandafter\let\csname c@#1\endcsname\c@dummy}

\cnewtheorem{theorem}[dummy]{Theorem}
\cnewtheorem{lemma}[dummy]{Lemma}
\cnewtheorem{corollary}[dummy]{Corollary}
\cnewtheorem{proposition}[dummy]{Proposition}

\cnewtheorem{conjecture}[dummy]{Conjecture}
\theoremstyle{definition}               %%Change Theoremstyle
\cnewtheorem{definition}[dummy]{Definition}
\cnewtheorem{question}[dummy]{Question}
\cnewtheorem{example}[dummy]{Example}
\newtheorem*{summaryun}{Summary}
\cnewtheorem{remark}[dummy]{Remark}
\newtheorem*{remarks}{Remarks}

\makeatother  %  move after \newtheorem block

\newcommand{\bz}{\mathbb Z}
\newcommand{\br}{\mathbb R}
\newcommand{\cc}{\mathcal C}
\newcommand{\cd}{\mathcal D}

\renewcommand{\to}{\rightarrow}
\let\barrsquare\square
\let\square\undefined

\begin{document}

\begin{abstract}
We prove a squeezing/stability theorem for delta-epsilon controlled
$L$--groups when the control map is a polyhedral stratified system of
fibrations on a finite polyhedron.  A relation with boundedly-controlled
$L$--groups is also discussed.
\end{abstract}

\maketitle

%%%%%%%%%%%%%%%%%%%%%%%%%%%%%%%%%%%%%%%%%%%%%%%%%%%%%%%%%%%%%%%%%%%%%%
\section{Introduction}
\label{section_intro}

Let us fix  an integer $n\ge 0$, a continuous map $p_X\co M\to X$ to a metric space $X$,
and a ring $R$ with involution.
For each pair of positive numbers $\epsilon\le\delta$,
the delta-epsilon controlled $L$--group
$L_n^{\delta,\epsilon}(X;p_X,R)$ is defined to be the set of equivalence classes of
$n$--dimensional quadratic Poincar\'e $R$--module complexes on $p_X$ of radius $\epsilon$
($=$ $n$--dimensional $\epsilon$--Poincar\'e
$\epsilon$--quadratic $R$--module
complexes on $p_X$),
where the equivalence relation is generated by Poincar\'e cobordisms of radius $\delta$
($=$ $\delta$--Poincar\'e $\delta$--cobordisms) -- see the work of Ranicki
and Yamasaki \cite{RY-L1,RY-L2,MY87}.
If $\delta\le\delta'$ and $\epsilon\le\epsilon'$, there is a natural homomorphism
\[
L_n^{\delta,\epsilon}(X;p_X,R)\to L_n^{\delta',\epsilon'}(X;p_X,R)
\]
defined by relaxation of control.
In general, this map is neither surjective nor injective.
None the less, if $X$ is a finite polyhedron and $p_X$ is a polyhedral
stratified system of fibrations in the sense of Quinn \cite{Q},
the map above turns out to be an isomorphism for certain values of $\delta$, $\delta'$,
$\epsilon$, $\epsilon'$:

\begin{theorem} [Stability in Controlled $L$--groups]
\label{stability} For each integer  $n\ge 0$ and a finite
polyhedron $X$, there exist constants $\delta_0>0$ and $\kappa>1$
such that the following hold: If
\begin{enumerate}
\item $\kappa\epsilon \le \delta \le \delta_0$,
$\kappa\epsilon' \le \delta' \le \delta'_0$,
$\delta\le\delta'$, $\epsilon\le\epsilon'$,
\item $p_X\co M\to X$ is a polyhedral stratified system of fibrations, and
\item $R$ is a ring with involution,
\end{enumerate}
then the relax-control map
$L_n^{\delta,\epsilon}(X;p_X,R)\to L_n^{\delta',\epsilon'}(X;p_X,R)$ is an isomorphism.
\end{theorem}

It follows that  all the groups $L_n^{\delta,\epsilon}(X;p_X,R)$ with
$\kappa\epsilon \le \delta\le \delta_0 $ are isomorphic and are equal to
the controlled $L$--group $L^c_n(X;p_X,R)$ of $p_X$ with coefficient ring $R$.

Stability is a consequence of squeezing;
squeezing/stability for controlled $K_0$ and $K_1$--groups were known (see
Pedersen \cite{EKP98}).
`Splitting' was the key idea there.
In \fullref{section_split}, we discuss splitting in the controlled $L$--theory.
An element of a controlled $L$--group is represented by a quadratic Poincar\'e
complex on a space.
If it splits into small pieces lying over cone-shaped sets ({\it e.g.} simplices),
then we can shrink all the pieces at the same time to obtain a squeezed complex.
But splitting in $L$--theory requires a change of $K$--theoretic decoration;
if you split a free quadratic Poincar\'e complex, then you get a projective
one in the middle.
Since the controlled reduced projective class group is known to vanish
when the coefficient ring is $\bz$ and the control map is $UV^1$,
we do not need to worry about the controlled $K$--theory and squeezing holds in this case
(see Pedersen--Quinn--Ranicki \cite{PQR}).

Several years ago the first named author proposed an approach to squeezing/stability
in controlled $L$--groups imitating the method of \cite{EKP98}.
The idea was to use projective complexes to split and to eventually
eliminate the projective pieces using the Eilenberg swindle:
\begin{align*}
[P]&=[P] + (-[P]+[P]) + (-[P]+[P]) + (-[P]+[P]) + \cdots \\
   &=([P]-[P]) + ([P]-[P]) + ([P]-[P]) + ([P]-[P]) + \cdots = 0 ~~.
\end{align*}
This approach works for any $R$ if $X$ is a circle;
we will briefly discuss the proof in \fullref{section_circle}.

The method used in \fullref{section_circle} does not generalize to higher dimensions,
because it requires repeated application of splitting but that is
not easy to do with projective complexes.
This means that we should not try to shrink the complex globally, but should
try to shrink a small part of the complex lying over a cone neighborhood of some point
at a time.
Such a local shrinking construction is possible when the control
map is a polyhedral stratified system of fibrations,
and is called an Alexander trick.  We study its effect in \fullref{section_alexander},
and  use it repeatedly to prove \fullref{stability} in \fullref{section_main_proof}.
Note that we do one splitting of the whole complex
for each application of an Alexander trick;
we are not splitting the split pieces.

In \fullref{section_variations}, we discuss several variations of \fullref{stability}.

Finally, in \fullref{boundedL}, we relate the delta-epsilon controlled
$L$--groups to the bounded $L$--groups in a special case.

The authors would like to thank Frank Connolly, Jim Davis, Frank Quinn and Andrew Ranicki
for invaluable suggestions.

%%%%%%%%%%%%%%%%%%%%%%%%%%%%%%%%%%%%%%%%%%%%%%%%%%%%%%%%%%%%%%%%%%%%%%
\section{Glueing and splitting}
\label{section_split}

In this section we review  techniques called glueing and splitting.
If $p_X\co M\to X$ is a control map and $Y$ is a subset of $X$, then we
denote the restriction $p_X|Y$ of $p_X$ by $p_Y$.  A closed $\epsilon$
neighborhood of $Y$ in $X$ is denoted by $Y^\epsilon$.  We refer the
reader to the papers \cite{RY-L1,RY-L2} by Ranicki and Yamasaki for
terms and notations in controlled $L$--theory.

We first discuss the glueing operation; it is to take the union of
two objects with common pieces of boundary.
Suppose there are consecutive Poincar\'e cobordisms of radius $\delta$,
one from $(C,\psi)$ to $(C',\psi')$ and the other from $(C',\psi')$ to $(C'',\psi'')$.
Then their union is a Poincar\'e cobordism of radius $100\delta$ from
$(C,\psi)$ to $(C'',\psi'')$ \cite[Proposition~2.8]{RY-L2}.
We will encounter this factor ``100'' many times in this article,
and will denote it by $\mu$ at several places of \fullref{section_main_proof}.
For example, we will need the following, which is a special case of
\cite[Proposition~3.7]{RY-L2}.

\begin{proposition}
\label{nullcobordism}
If $[C,\psi]=0$ in $L_n^{\delta,\epsilon}(X;p_X,R)$,
then there is a Poincar\'e cobordism of radius $100\delta$ from $(C,\psi)$ to $0$.
\end{proposition}

\begin{proof}
By definition, there is a sequence of consecutive Poincar\'e cobordisms
starting from $(C,\psi)$ and ending at $0$. Their union can be regarded
as the union of the even-numbered ones and the odd-numbered ones, so it
is $100\delta$ Poincar\'e.
\end{proof}

Next we discuss splitting.
Before stating the splitting lemma, let us recall a minor technicality from
Ranicki--Yamasaki \cite[Section~6]{RY-K}: Suppose $X$ is the union of two closed subsets $A$ and $B$ with
intersection $Y=A\cap B$.  If a path $\gamma\co [0,s]\to M$ with $p_X\gamma(0)\in A$
is contained in $p_X^{-1}(\{\gamma(0)\}^\epsilon)$, then it lies in
$p_X^{-1}(A\cup Y^{2\epsilon})$.
Of course it is contained also in $p_X^{-1}(A^{\epsilon})$, but this is slightly
less useful.

\begin{lemma}[Splitting Lemma]
\label{splitting_lemma}
For any integer $n\ge 2$,
there exists a positive number $\lambda\ge 1$ such that the following holds:
If $p_X\co M\to X$ is a map to a metric space $X$,
$X$ is the union of two closed subsets $A$ and $B$ with intersection $Y$,
and $R$ is a ring with involution, then
for any $n$--dimensional quadratic Poincar\'e $R$--module complex
$c=(C,\psi)$ on $p_X$ of radius $\epsilon$,
there exist a Poincar\'e cobordism of radius $\lambda\epsilon$ from $c$ to the union
$c'\cup c''$ of an $n$--dimensional quadratic Poincar\'e pair
$c'=(f'\co P\to C', (\delta\wbar\psi', -\wbar\psi))$
on $p_{A\cup Y^{\lambda\epsilon}}$ of radius $\lambda\epsilon$ and
an $n$--dimensional quadratic Poincar\'e pair
$c''=(f''\co P\to C'', (\delta\wbar\psi'', \wbar\psi))$
on $p_{B\cup Y^{\lambda\epsilon}}$ of radius $\lambda\epsilon$,
where $(P,\wbar\psi)$ is an $(n-1)$--dimensional quadratic Poincar\'e projective
$R$--module complex on $p_{Y^{\lambda\epsilon}}$
and $P$ is $\lambda\epsilon$ chain equivalent to an
$(n-1)$--dimensional free chain complex on $p_{A\cup Y^{\lambda\epsilon}}$
and also to an
$(n-1)$--dimensional free chain complex on $p_{B\cup Y^{\lambda\epsilon}}$.
\end{lemma}

\begin{proof} This is an epsilon-control version of Ranicki's
argument for the bounded control case \cite{LowerKL}.
For a given $(C,\psi)$ of radius $\epsilon$, pick up a subcomplex
$C'\subset C$ such that $C'$ is identical with $C$ over $A$ and $C'$ lies over
some neighborhood of $A$.
Let $p\co C\to C/C'$ be the quotient map and define $C''$ by the $n$--dual $(C/C')^{n-*}$.
Define a complex $E$ by the desuspension $\Omega\cc(p\cd_\psi p^*)$
of the algebraic mapping cone of the map
\[
\xymatrix@1{
C''=(C/C')^{n-*} \ar[r]^-{p^*} & C^{n-*} \ar[r]^-{\cd_\psi} & C \ar[r]^-{p} & C/C'},
\]
where $\cd_\psi$ is the duality map $(1+T)\psi_0$ for $\psi$.
There are natural maps $g'\co E\to C'$, $g''\co E\to C''$ and adjoining
quadratic Poincar\'e structures on them such that the union along
the common boundary is homotopy equivalent to the original complex
$c$. We should note that $E$ is non-trivial in degrees $-1$ and
$n$ and that it lies over $B$.

Since $\cd_\psi$ is a small chain equivalence, its mapping cone is
contractible. Therefore, $E$ is contractible away from the union
of $A$ and a small neighborhood of $Y$, and it is chain equivalent
to a projective chain complex $P$ lying over a small neighborhood
of $Y$ by \cite[Sections~5.1 and 5.2]{RY-K}. Note that $\bz$ is used as
the coefficient ring in \cite{RY-K}, but the same argument works
when the coefficient ring is replaced by $R$. Since $n\ge 2$, we
can assume that $P$ is strictly $(n-1)$--dimensional ({\it i.e.}
$C_i=0$ for $i<0$ and $i>n-1$ ) by the standard folding argument,
and the chain equivalence induces a desired cobordism.

There is a quadratic Poincar\'e structure on a chain map $f'\co P\to C'$;
therefore, the duality map gives a chain equivalence
$C'^{n-*} \longrightarrow {\mathcal C}(f')$,
were ${\mathcal C}(f')$ denotes the algebraic mapping cone of $f'\co P\to C'$.
Therefore
\[
[P]=-([C'] - [P]) = -[{\mathcal C}(f')] = -[C'^{n-*}] = 0
\]
in the epsilon controlled reduced projective class group of the union of $A$ and
a small neighborhood of $Y$ with coefficient in $R$,
and hence $P$ is chain equivalent to a free $(n-1)$--dimensional complex $F'$ lying
over the union of $A$ and a small neighborhood of $Y$.
\end{proof}

\begin{remarks}
(1)\qua
Suppose that $X$ is a finite polyhedron or a finite cell complex in the
sense of Rourke and Sanderson \cite{RS72}
more generally.
Then there exist positive numbers $\epsilon_X>0$, $\mu_X\ge 1$ and a homotopy
$\{f_t\}\co X\to X$ such that
\begin{itemize}
\item $f_0=1_X$,
\item $f_t(\Delta)\subset \Delta$ for each cell $\Delta$ and for each $t\in[0,1]$,
\item $f_t$ is Lipschitz with Lipschitz constant $\mu_X$ for each $t\in[0,1]$, and
\item $f_1((X^{(i)})^{\epsilon_X})\subset X^{(i)}$ for every $i$, where
$X^{(i)}$ is the $i$--skeleton of $X$.
\end{itemize}
Suppose $\{f_t\}$ is covered by a homotopy $\{F_t\co M\to M\}$ and
set $\delta_X=\epsilon_X\lambda$.  If $\epsilon\le\delta_X$ and $A$ and $B$ are subcomplexes
of $X$, then by applying $F_1$ to the splitting given in the above lemma,
we may assume that the pieces lie over $A$, $B$ and  $A\cap B$ respectively,
instead of their neighborhoods, since the homotopy gives small isomorphisms
between the corresponding pieces.
But $\lambda$ is now replaced by $\mu_X\lambda$ and it depends not only on $n$ but also on $X$.
We call such a deformation $\{f_t\}$ a {\it rectification} for $X$.

(2)\qua The splitting formula for pairs given by Yamasaki \cite{MY87}
can be combined with
\cite[Sections~5.1 and~5.2]{RY-K} to prove a similar splitting lemma for pairs of dimension $\ge 3$:
a sufficiently small Poincar\'e pair splits into two adjoining quadratic Poincar\'e triads
whose common boundary piece is possibly a projective pair.
\end{remarks}

%%%%%%%%%%%%%%%%%%%%%%%%%%%%%%%%%%%%%%%%%%%%%%%%%%%%%%%%%%%%%%%%%%%%%%
\section{Squeezing over a circle}
\label{section_circle}
We discuss squeezing over the unit circle.
We use the maximum metric of $\br^2$, so the unit circle looks like a square:

\centerline{\includegraphics{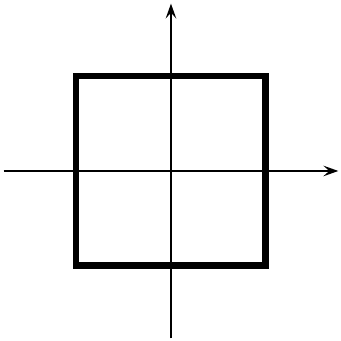}}

Consider a quadratic Poincar\'e $R$--module complex on the unit circle.
We assume that its radius is sufficiently small so that it splits
into four free pieces $E$, $F$, $G$, $H$ with projective boundary
pieces $P$, $Q$, $S$, $T$ as shown in the picture below.
The shadowed region is a cobordism between the original complex and the union of
$E$, $F$, $G$, $H$.
Although we actually measure the radius using the radial projection to the unit circle
({i.e.} the square), we pretend that complexes and cobordisms are over the plane.

{\labellist\small
\hair=4pt
\pinlabel {$T$} [br] at 92 700
\pinlabel {$E$} [bl] at 134 700
\pinlabel {$P$} [bl] at 176 700
\pinlabel {$H$} [br] at 92 658
\pinlabel {$F$} [bl] at 176 658
\pinlabel {$S$} [tr] at 92 616
\pinlabel {$G$} [tl] at 134 616
\pinlabel {$Q$} [tl] at 176 616
\endlabellist
\centerline{\includegraphics{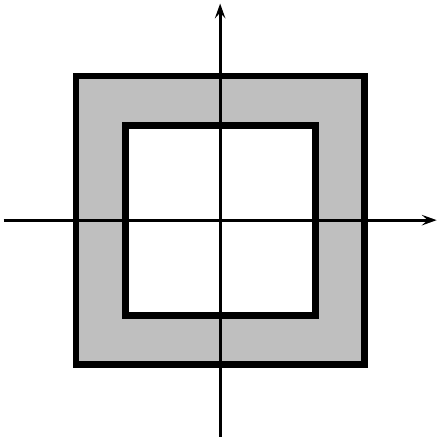}}}

We extend this cobordism in the following way. On the right
vertical edge, we have a quadratic pair $P\oplus Q\to F$. (We are
omitting the quadratic structure from notation.) Take the tensor
product of this with the symmetric complex of the unit interval
$[0,1]$. Make many copies of such a product and consecutively glue
them one after the other to the cobordism. Do the same thing with
the other three edges. Then fill in the four quadrants by copies
of $P$, $Q$, $S$, $T$ multiplied by the symmetric complex of
$[0,1]^2$ so that the whole picture looks like a huge square with
a square hole at the center.

Although this cobordism is made up of free complexes and projective complexes,
the projective complexes sitting on the white edges are shifted up 1 dimension,
and the projective complexes sitting at the lattice points are
shifted up 2 dimension in the union.
$$
\labellist
\tiny
\hair=2pt

\pinlabel {$T$} [br] at  88 704
\pinlabel {$T$} [r]  at  88 693
\pinlabel {$T$} [r]  at  88 682
\pinlabel {$T$} [r]  at  88 671
\pinlabel {$T$} [r]  at  88 660
\pinlabel {$T$} [r]  at  88 648
\pinlabel {$T$} [r]  at  88 636

\pinlabel {$T$} [b]  at 100 704
\pinlabel {$T$}      at 100 693
\pinlabel {$T$}      at 100 682
\pinlabel {$T$}      at 100 671
\pinlabel {$T$}      at 100 660
\pinlabel {$T$}      at 100 648
\pinlabel {$T$}      at 100 636

\pinlabel {$T$} [b]  at 110 704
\pinlabel {$T$}      at 110 693
\pinlabel {$T$}      at 110 682
\pinlabel {$T$}      at 110 671
\pinlabel {$T$}      at 110 660
\pinlabel {$T$}      at 110 648
\pinlabel {$T$}      at 110 636

\pinlabel {$T$} [b]  at 122 704
\pinlabel {$T$}      at 122 693
\pinlabel {$T$}      at 122 682
\pinlabel {$T$}      at 122 671
\pinlabel {$T$}      at 122 660
\pinlabel {$T$}      at 122 648
\pinlabel {$T$}      at 122 636

\pinlabel {$T$} [b]  at 133 704
\pinlabel {$T$}      at 133 693
\pinlabel {$T$}      at 133 682
\pinlabel {$T$}      at 133 671
\pinlabel {$T$}      at 133 660
\pinlabel {$T$}      at 133 648
\pinlabel {$T$}      at 133 636

\pinlabel {$T$} [b]  at 144 704
\pinlabel {$T$}      at 144 693
\pinlabel {$T$}      at 144 682
\pinlabel {$T$}      at 144 671
\pinlabel {$T$}      at 144 660
\pinlabel {$T$}      at 144 648
\pinlabel {$T$}      at 144 636

\pinlabel {$T$} [b]  at 156 704
\pinlabel {$T$}      at 156 693
\pinlabel {$T$}      at 156 682
\pinlabel {$T$}      at 156 671
\pinlabel {$T$}      at 156 660
\pinlabel {$T$}      at 156 648
\pinlabel {$T$}      at 156 636

\pinlabel {$E$} [b]  at 167 704
\pinlabel {$E$}      at 167 693
\pinlabel {$E$}      at 167 682
\pinlabel {$E$}      at 167 671
\pinlabel {$E$}      at 167 660
\pinlabel {$E$}      at 167 648
\pinlabel {$E$}      at 167 636

\pinlabel {$P$} [b]  at 179 704
\pinlabel {$P$}      at 179 693
\pinlabel {$P$}      at 179 682
\pinlabel {$P$}      at 179 671
\pinlabel {$P$}      at 179 660
\pinlabel {$P$}      at 179 648
\pinlabel {$P$}      at 179 636

\pinlabel {$P$} [b]  at 190 704
\pinlabel {$P$}      at 190 693
\pinlabel {$P$}      at 190 682
\pinlabel {$P$}      at 190 671
\pinlabel {$P$}      at 190 660
\pinlabel {$P$}      at 190 648
\pinlabel {$P$}      at 190 636

\pinlabel {$P$} [b]  at 201 704
\pinlabel {$P$}      at 201 693
\pinlabel {$P$}      at 201 682
\pinlabel {$P$}      at 201 671
\pinlabel {$P$}      at 201 660
\pinlabel {$P$}      at 201 648
\pinlabel {$P$}      at 201 636

\pinlabel {$P$} [b]  at 212 704
\pinlabel {$P$}      at 212 693
\pinlabel {$P$}      at 212 682
\pinlabel {$P$}      at 212 671
\pinlabel {$P$}      at 212 660
\pinlabel {$P$}      at 212 648
\pinlabel {$P$}      at 212 636

\pinlabel {$P$} [b]  at 224 704
\pinlabel {$P$}      at 224 693
\pinlabel {$P$}      at 224 682
\pinlabel {$P$}      at 224 671
\pinlabel {$P$}      at 224 660
\pinlabel {$P$}      at 224 648
\pinlabel {$P$}      at 224 636

\pinlabel {$P$} [b]  at 235 704
\pinlabel {$P$}      at 235 693
\pinlabel {$P$}      at 235 682
\pinlabel {$P$}      at 235 671
\pinlabel {$P$}      at 235 660
\pinlabel {$P$}      at 235 648
\pinlabel {$P$}      at 235 636

\pinlabel {$P$} [l]  at 246 704
\pinlabel {$P$} [l]  at 246 693
\pinlabel {$P$} [l]  at 246 682
\pinlabel {$P$} [l]  at 246 671
\pinlabel {$P$} [l]  at 246 660
\pinlabel {$P$} [l]  at 246 648
\pinlabel {$P$} [l]  at 246 636

\pinlabel {$H$} [r]  at  88 625
\pinlabel {$H$}      at 100 625
\pinlabel {$H$}      at 110 625
\pinlabel {$H$}      at 122 625
\pinlabel {$H$}      at 133 625
\pinlabel {$H$}      at 144 625
\pinlabel {$H$}      at 156 625

\pinlabel {$F$}      at 179 625
\pinlabel {$F$}      at 190 625
\pinlabel {$F$}      at 201 625
\pinlabel {$F$}      at 212 625
\pinlabel {$F$}      at 224 625
\pinlabel {$F$}      at 235 625
\pinlabel {$F$} [l]  at 246 625

\pinlabel {$S$} [r]  at  88 614
\pinlabel {$S$} [r]  at  88 603
\pinlabel {$S$} [r]  at  88 591
\pinlabel {$S$} [r]  at  88 580
\pinlabel {$S$} [r]  at  88 569
\pinlabel {$S$} [r]  at  88 556
\pinlabel {$S$} [tr] at  88 546

\pinlabel {$S$}      at 100 614
\pinlabel {$S$}      at 100 603
\pinlabel {$S$}      at 100 591
\pinlabel {$S$}      at 100 580
\pinlabel {$S$}      at 100 569
\pinlabel {$S$}      at 100 556
\pinlabel {$S$} [t]  at 100 546

\pinlabel {$S$}      at 110 614
\pinlabel {$S$}      at 110 603
\pinlabel {$S$}      at 110 591
\pinlabel {$S$}      at 110 580
\pinlabel {$S$}      at 110 569
\pinlabel {$S$}      at 110 556
\pinlabel {$S$} [t]  at 110 546

\pinlabel {$S$}      at 122 614
\pinlabel {$S$}      at 122 603
\pinlabel {$S$}      at 122 591
\pinlabel {$S$}      at 122 580
\pinlabel {$S$}      at 122 569
\pinlabel {$S$}      at 122 556
\pinlabel {$S$} [t]  at 122 546

\pinlabel {$S$}      at 133 614
\pinlabel {$S$}      at 133 603
\pinlabel {$S$}      at 133 591
\pinlabel {$S$}      at 133 580
\pinlabel {$S$}      at 133 569
\pinlabel {$S$}      at 133 556
\pinlabel {$S$} [t]  at 133 546

\pinlabel {$S$}      at 144 614
\pinlabel {$S$}      at 144 603
\pinlabel {$S$}      at 144 591
\pinlabel {$S$}      at 144 580
\pinlabel {$S$}      at 144 569
\pinlabel {$S$}      at 144 556
\pinlabel {$S$} [t]  at 144 546

\pinlabel {$S$}      at 156 614
\pinlabel {$S$}      at 156 603
\pinlabel {$S$}      at 156 591
\pinlabel {$S$}      at 156 580
\pinlabel {$S$}      at 156 569
\pinlabel {$S$}      at 156 556
\pinlabel {$S$} [t]  at 156 546

\pinlabel {$G$}      at 167 614
\pinlabel {$G$}      at 167 603
\pinlabel {$G$}      at 167 591
\pinlabel {$G$}      at 167 580
\pinlabel {$G$}      at 167 569
\pinlabel {$G$}      at 167 556
\pinlabel {$G$} [t]  at 167 546

\pinlabel {$Q$}      at 179 614
\pinlabel {$Q$}      at 179 603
\pinlabel {$Q$}      at 179 591
\pinlabel {$Q$}      at 179 580
\pinlabel {$Q$}      at 179 569
\pinlabel {$Q$}      at 179 556
\pinlabel {$Q$} [t]  at 179 546

\pinlabel {$Q$}      at 190 614
\pinlabel {$Q$}      at 190 603
\pinlabel {$Q$}      at 190 591
\pinlabel {$Q$}      at 190 580
\pinlabel {$Q$}      at 190 569
\pinlabel {$Q$}      at 190 556
\pinlabel {$Q$} [t]  at 190 546

\pinlabel {$Q$}      at 201 614
\pinlabel {$Q$}      at 201 603
\pinlabel {$Q$}      at 201 591
\pinlabel {$Q$}      at 201 580
\pinlabel {$Q$}      at 201 569
\pinlabel {$Q$}      at 201 556
\pinlabel {$Q$} [t]  at 201 546

\pinlabel {$Q$}      at 212 614
\pinlabel {$Q$}      at 212 603
\pinlabel {$Q$}      at 212 591
\pinlabel {$Q$}      at 212 580
\pinlabel {$Q$}      at 212 569
\pinlabel {$Q$}      at 212 556
\pinlabel {$Q$} [t]  at 212 546

\pinlabel {$Q$}      at 224 614
\pinlabel {$Q$}      at 224 603
\pinlabel {$Q$}      at 224 591
\pinlabel {$Q$}      at 224 580
\pinlabel {$Q$}      at 224 569
\pinlabel {$Q$}      at 224 556
\pinlabel {$Q$} [t]  at 224 546

\pinlabel {$Q$}      at 235 614
\pinlabel {$Q$}      at 235 603
\pinlabel {$Q$}      at 235 591
\pinlabel {$Q$}      at 235 580
\pinlabel {$Q$}      at 235 569
\pinlabel {$Q$}      at 235 556
\pinlabel {$Q$} [t]  at 235 546

\pinlabel {$Q$} [l]  at 246 614
\pinlabel {$Q$} [l]  at 246 603
\pinlabel {$Q$} [l]  at 246 591
\pinlabel {$Q$} [l]  at 246 580
\pinlabel {$Q$} [l]  at 246 569
\pinlabel {$Q$} [l]  at 246 556
\pinlabel {$Q$} [tl] at 246 546
\endlabellist
\centerline{\includegraphics{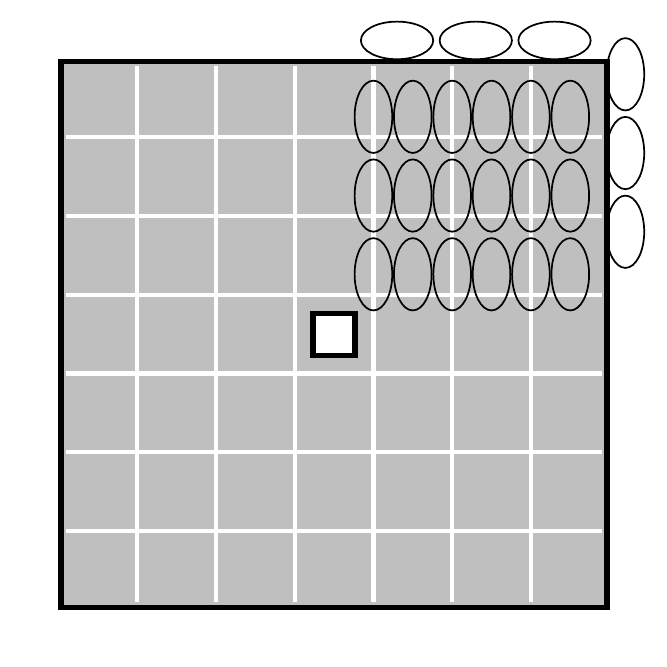}}
$$
We can make pairs of these (as shown in the picture above for $P$'s)
so that each pair contributes the trivial element in the controlled
reduced projective class group.
Replace each pair by a free complex.

Unlike the real Eilenberg swindle, there are four projective complexes
left which do not make pairs.
We may assume that they are the boundary pieces of $F$ and $H$ on the outer end.
Since the two pairs $P\oplus Q\to F$, $S\oplus T\to H$ are Poincar\'e,
the unions $P\oplus Q$ and $S\oplus T$ are locally chain equivalent to free complexes.
Thus we can replace them by free complexes, and now everything is free.

Now recall that we actually measure things by a radial projection to the square.
Thus we have a cobordism from the original complex to another complex of very small radius.
If we increase the number of layers in the construction, the radius of the outer end
becomes arbitrarily small.
This is the squeezing in the case of $S^1$.

%%%%%%%%%%%%%%%%%%%%%%%%%%%%%%%%%%%%%%%%%%%%%%%%%%%%%%%%%%%%%%%%%%%%%%
\section{Alexander trick}
\label{section_alexander}

The method in the previous section does not work for higher dimensional complexes,
because we cannot inductively split the projective pieces.
But the proof suggests an alternative way toward squeezing/stability.
This is the topic of this section.
Although we used a radial projection to measure the size in the previous section,
we draw pictures of things in their real sizes in this section.

Let us fix an integer $n\ge 2$ and a finite polyhedron $X$.
All the complexes below are $R$--module complexes, where $R$ is
a ring with involution.
We assume that the control map $p_X\co M\to X$ is a polyhedral stratified
system of fibrations (see Quinn \cite{Q});
$p_X$ is fiber homotopy equivalent to a map $q_X\co N\to X$ which has an iterated
mapping cylinder decomposition in the sense of Hatcher \cite{Hatcher}:
there is a partial order on the set of the vertices of $X$ such that,
for each simplex $\Delta$ of $X$,
\begin{enumerate}
\item the partial order restricts to a
total order of the vertices of $\Delta$
\[
v_0<v_1<\dots<v_k~,
\]
\item $q_X^{-1}(\Delta)$ is the iterated mapping cylinder of a sequence of maps
\[
F_{v_0}\longrightarrow F_{v_1}\longrightarrow \dots\longrightarrow F_{v_k}~,
\]
\item the restriction $q_X|q_X^{-1}(\Delta)$ is the natural map induced from
the iterated mapping cylinder structure of $q_X^{-1}(\Delta)$ above and the
iterated mapping cylinder structure of $\Delta$ coming from the sequence
\[
\{v_0\}\longrightarrow \{v_1\}\longrightarrow \dots \longrightarrow\{v_k\}~.
\]
\end{enumerate}
An order on the set of the vertices of $X$ is said to be \emph{compatible with}
$p_X$ if it is compatible with this partial order.  Let us fix an order compatible
with $p_X$.

Pick a vertex $v$ of $X$, and let $A$ be the star neighborhood of $v$,
$B$ be the closure of the complement of $A$ in $X$,
and $S$ be the union of the simplices in $A$ whose vertices are all $\ge v$
with respect to the chosen order.
This will be called the {\it stable set} at $v$.
Let $s\co A\to S$ be the simplicial retraction defined by
\[
s(v')=\begin{cases}
v &\hbox{\qquad if $v'< v$,}\\
v' &\hbox{\qquad if $v'\ge v$,}
\end{cases}
\]
for vertices $v'$ of $A$.  A strong deformation retraction $s_t\co A\to A$ is defined
by $s_t(a)=(1-t)a+t\,s(a)$ for $a\in A$ and $t\in [0,1]$.
Note that this strong deformation retraction $s_t$ is covered by a deformation
$\wtilde s_t$ on $M$, since $p_X$ is a polyhedral stratified system of fibrations.

{\small
\labellist
\pinlabel {$1$} [r]  at 87 636
\pinlabel {$2$} [tl] at 210 601
\pinlabel {$3$} [br] at 155 636
\pinlabel {$4$} [b]  at 140 706
\pinlabel {$5$} [t]  at 156 581
\pinlabel {$6$} [l]  at 228 669
\pinlabel {$S$}      at 174 666
\pinlabel {$S$} [r]  at 157 610
\endlabellist
\centerline{\includegraphics{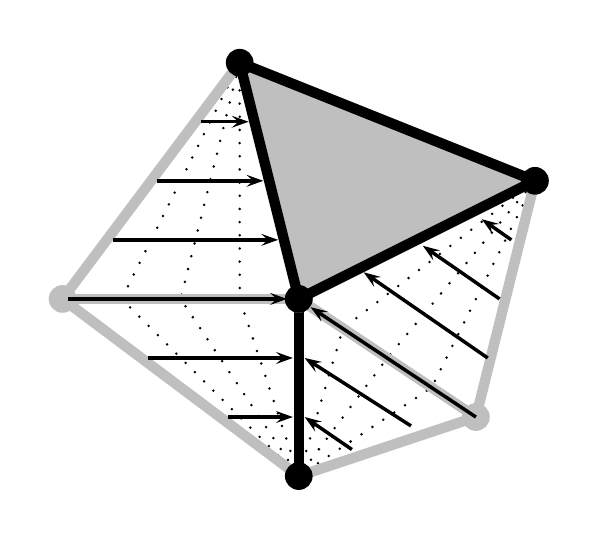}}}

Given a sufficiently small $n$--dimensional quadratic Poincar\'e complex
$c=(C,\psi)$ on $p_X$, one can split it according to the splitting of $X$
into $A$ and $B$: $c$ is cobordant (actually homotopy equivalent)
to the union $c'$ of a projective quadratic Poincar\'e pair
$a=(f\co P\to F, (\delta\psi', \psi'))$ on $p_A$ and a projective quadratic Poincar\'e pair
$b=(g\co P\to G, (\delta\psi'',-\psi'))$ on $p_B$, where
$F$ is an $n$--dimensional chain complex on $p_A$,
$G$ is an $n$--dimensional chain complex on $p_B$,
and $P$ is an $(n-1)$--dimensional projective chain complex on $p_{A\cap B}$.
Here we again used the assumption on $p_X$.
See the remark after the splitting lemma.

Make many copies of the product cobordism from the pair $a$ to itself,
and successively glue them to the cobordism between $c$ and $c'$.
This gives us a cobordism from $c$ to a (possibly) projective complex
as in the left picture below.

\begin{center}
\includegraphics{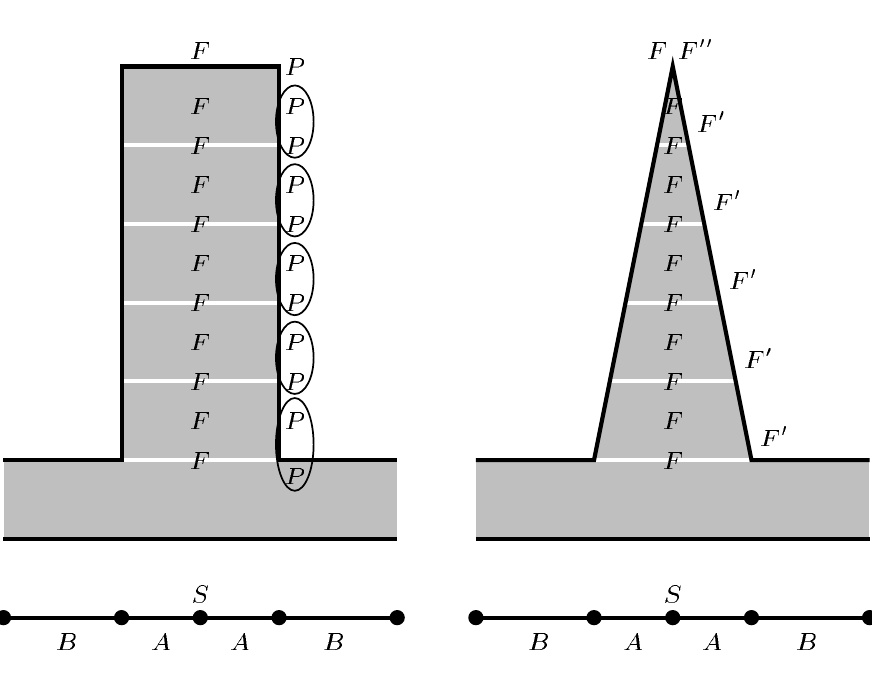}
\end{center}

We will remedy the situation by replacing the projective end by
a free complex as follows.
The copies of $P$ connecting the layers are actually shifted up 1 dimension
in the union,
so the marked pairs of $P$'s contribute the trivial element of the controlled
$\wtilde K_0$ group of $A\cap B$, and
we can replace each pair with a free module by adding chain complexes
of the form
\[
\xymatrix@1{
0\ar[r] & Q_i \ar[r]^{1} &Q_i \ar[r] & 0}
\]
lying over  $A\cap B$,
where $Q_i$ is a projective module such that $P_i\oplus Q_i$ is free.
Therefore, these pairs are all chain equivalent to some free chain  complex $F'$.
The last $P$ remaining at the top of the picture can be replaced
by some free complex $F''$ lying over $A$ as stated in the splitting lemma.

We deform the tower, which is now free, toward $S$ using $\wtilde s_t$ as in the picture
above so that the top of the tower is  completely deformed to $S$.

\begin{summaryun}  {\it There exist constants $\delta>0$ and $\lambda\ge 1$
which depend on $n$ and $X$
such that any $n$--dimensional quadratic Poincar\'e complex of
radius $\epsilon\le\delta$ on $p_X$ is
$\lambda\epsilon$ Poincar\'e cobordant to another complex which is
small in the track direction of $s_t$.
The more layers we use, the smaller the result becomes in the track direction.
}
\end{summaryun}

\begin{remarks}
(1)\qua We cannot take $\lambda=1$ in general,
since the radius of the complexes gets bigger during the splitting/glueing processes.

(2)\qua This construction will be referred to as the {\it Alexander trick} at $v$.

(3)\qua There is also an {\it Alexander trick for pairs}.
If we use the splitting lemma for pairs, then instead of a pair we get a Poincar\'e triad
\[
\xymatrix@1{
P \ar[r] \ar[d] \ar@{~>}[rd] & Q \ar[d] \\
E \ar[r] & F}
\]
over $A$, where $P$, $Q$ are projective and $E$, $F$ are free.
Since both $P$ and $Q$ are free over $A$, we can carry out the construction
exactly in the same manner as above.
The effect on the boundary is exactly the same as the absolute Alexander trick.

(4)\qua Take a simplex $\Delta$ of $X$ with ordered vertices
$v_0<v_1<\cdots < v_n$. Let $(\lambda_0, \dots, \lambda_n)$ be the
barycentric coordinates of a point $x\in\Delta$, {\it i.e.}
$x=\sum \lambda_i v_i$. Then we define the {\it
pseudo-coordinates} $(x_1,\dots,x_n)$ of $x$ by
$x_i=\lambda_i/(\lambda_0+\cdots + \lambda_i)$. Actually $x_i$ is
indeterminate if $\lambda_0=\cdots=\lambda_i$. Let
$s_{i,t}\co \Delta\to \Delta$ be the restriction to $\Delta$ of the
deformation retraction used for an Alexander trick at $v_i$; then
$s_{0,t}=1_\Delta$ for every $t\in[0,1]$, and $s_{i,t}$ preserves
the pseudo-coordinate $x_j$ for $j$ not equal to $i$. This means
that, roughly speaking, an Alexander trick at $v_i$ improves the
radius control in the $x_i$ direction and changes the radius
control in the $x_j$ direction ($j\ne i$) only up to
multiplications by the constant $\lambda$ given in the Splitting
Lemma and by the Lipschitz constant of $s_{i,t}$ which is uniform
with respect to $t$. Thus, if we can perform appropriate Alexander
tricks at all the vertices of $\Delta$, then we can obtain an
arbitrarily fine control over $\Delta$. A more detailed discussion
will be given in the next section.
\end{remarks}

Let us state a lemma on Lipschitz properties related to the homotopy $s_t$ above,
for future use.

\begin{lemma}
\label{lipschitz}
Let $X$ be a subset of $\br^N$ with diameter $d$ and $s\co X\to X$ be a Lipschitz map
with Lipschitz constant $K\ge 1$.  Suppose $X$ contains the line segment $xs(x)$
for every $x\in X$ and let $s_t(x)=ts(x)+(1-t)x$ for $t\in[0,a]$.
Then $s_t\co X\to X$ has Lipschitz constant $K$, and
the map
\[
H\co X\times[0,a]\to X\times[0,a]~;\quad H(x,t)=(s_{t/a}(x), t)
\]
has Lipschitz constant $\max\{d/a, 1\}+K$ with respect to the maximum metric on
$X\times[0,a]$.
\end{lemma}

\begin{proof} Let $x$, $y$ be points in $X$.  Then
\begin{align*}
d(s_t(x),s_t(y)) &= \lVert t(s(x)-s(y))+(1-t)(x-y)\rVert\\
&\le t d(s(x),s(y)) + (1-t)d(x,y)\\
&\le tK d(x,y) + (1-t)Kd(x,y)=Kd(x,y)~.
\end{align*}
Next, take two points $p=(x,t)$, $q=(y,u)$ of $X\times[0,a]$, and
let $p'=(x,u)$.
Then we have
\begin{align*}
d(H(p),H(q))&\le d(H(p),H(p'))+d(H(p'),H(q))\\
&= \max\{d(s_{t/a}(x),s_{u/a}(x)), |t-u|\} + d(s_{u/a}(x),s_{u/a}(y))\\
&\le  \max\{|t-u|d(s(x),x)/a, |t-u|\} + K d(x,y)\\
&\le  |t-u|\max\{d/a, 1\} + K d(x,y)\\
&\le  (\max\{d/a, 1\} + K)\max\{d(x,y), |t-u|\}\\
&=(\max\{d/a, 1\} + K)d(p,q)
\end{align*}
which completes the proof.
\end{proof}

%%%%%%%%%%%%%%%%%%%%%%%%%%%%%%%%%%%%%%%%%%%%%%%%%%%%%%%%%%%%%%%%%%%%%%
\section[Proof of Theorem \ref{stability}]{Proof of \fullref{stability}}
\label{section_main_proof}

The algebraic theory of surgery on quadratic Poincar\'e complexes in an additive
category (see Ranicki \cite{AdditiveL}) carries over nicely to the
controlled setting, and can be used to prove a stable periodicity of the
controlled $L$--groups.  Therefore, we give a proof of the stability in
the case $n\ge 2$.  The stability for $n=0, 1$ follows from the stability
for $n=4, 5$.

We first state the squeezing lemma for quadratic Poincar\'e complexes:

\begin{lemma}[Squeezing of Quadratic Poincar\'e Complexes]
\label{squeeze}
Let $n\ge 2$ be an integer and $X$ be a finite polyhedron.
There exist constants $\delta_0>0$ and $\kappa>1$
such that the following hold: If $\epsilon<\epsilon'\le\delta_0$, then
any $n$--dimensional quadratic Poincar\'e $R$--module complex of
radius $\epsilon'$ on a polyhedral stratified system of fibrations
over $X$ is $\kappa\epsilon'$ Poincar\'e cobordant
to a quadratic Poincar\'e complex of radius $\epsilon$.
\end{lemma}

%%%%%% main proof %%%%%%%
\begin{proof}
Let $X$ be a polyhedron in $\br^N$, and
$p_X\co M\to X$ be a polyhedral stratified system of fibrations.
Order the vertices of $X$ compatibly with $p_X$:
\[
v_0 < v_1 < \dots < v_m
\]
The basic idea is to apply the Alexander trick at each $v_i$.
This should make the complex arbitrarily small in $X$ as noted in the previous section.
The problem is that an Alexander trick is made up of two steps:
the first step is to make a tower using splitting, and the second step
is to squeeze the tower, and estimating the effect of the splitting used in the first step
is very difficult especially near the vertex when the object is getting smaller in a non-uniform way.
To avoid this difficulty, we blow up the metric around each vertex so that
the ordinary control on the new metric space insures us that the result
has a desired small control measured on the original metric space $X$.
Note that we are implicitly using this approach in the circle case.

Let us start from a complex $c$ of radius $\epsilon'>0$ on $X$. Since $X$
is a finite polyhedron, there exist $\delta>0$ and $\lambda\ge1$
such that if $\epsilon'\le\delta$ then $c$ is $\lambda\epsilon'$
cobordant to the union of two pieces according to the splitting of
$X$ into two subpolyhedra as in the remark after \fullref{splitting_lemma}. Recall that $\delta$ and $\lambda$ depends
on $X$. Set $\mu=100$, and set $\delta_0=\delta /
(\mu\lambda^2)^{m-1}$. The factor 100 comes from
\cite[Section~2.8]{RY-L2} as was mentioned in \fullref{section_split}. 
We claim that if $\epsilon'\le\delta_0$, then a
successive application of Alexander tricks produces a cobordism
from $c$  to a complex of radius $\epsilon$.

Let us fix some more notation.
$V_1$, \dots, $V_m$ are the star neighborhoods of $v_1$, \dots, $v_m$,
and $L_1$, \dots, $L_m$ are the links of $v_1$, \dots, $v_m$;
$V_i$ is the cone over $L_i$ with vertex $v_i$ for each $i$.
$S_1$, \dots, $S_m$ are the stable sets at $v_1$, \dots, $v_m$.
$K \ge 1$ is the Lipschitz constant which works for every retraction
$s_{i}\co V_i\to S_i$ used for the Alexander trick at $v_i$.
Let $d$ denote the diameter of $X$, and let
$\sharp(X)$ denote the number of simplices of $X$.
Now fix a number $H\ge 1$ such that
\[
H > d ~ \quad \text{and} \quad
4\mu\sharp(X)(K+1)^m(\mu\lambda^2)^m\epsilon'\cdot \displaystyle \frac{d}{~H~} <\epsilon~.
\]
We inductively define metric spaces and subsets
\begin{gather*}
 X^{i,j}_* \supset  X^{i,j}
\supset  V^{i,j}_{k}\supset  L^{i,j}_{k}
\qquad (1\le i\le j<k\le m)
\end{gather*}
together with control maps
$ p^{i,j}_*\co  M^{i,j}_*\to  X^{i,j}_*$
as follows.

Identify $\br^N$ with the subset $\br^N\times\{0\}$ of $\br^{N+1}=\br^N\times\br$
with the maximum product metric.
For each $i=1,\dots,m$,
define  $X^{i,i}_*$ and its subsets $X^{i,i}$, $V^{i,i}_k$, $L^{i,i}_k$ ($k=i+1,\cdots,m$) by
\begin{align*}
 X^{i,i}_*&=X \cup (V_i\times [0,H]),\\
 X^{i,i} &=(X-V_i)\cup (L_i\times [0,H]) \cup (V_i\times \{H\}),\\
 V^{i,i}_k &=(V_k-V_i)\cup (V_k\cap L_i\times [0,H]) \cup (V_k\cap V_i\times \{H\}),\\
 L^{i,i}_k &=(L_k-V_i)\cup (L_k\cap L_i\times [0,H]) \cup (L_k\cap V_i\times \{H\}).
\end{align*}
The projection $\br^N\times\br\to \br^N$ restricts to a retraction
$r_{i,i}\co X^{i,i}_*\to X$.
We define the control map
$ p^{i,i}_* \co   M^{i,i}_*\to X^{i,i}_*$
to be the pullback of $p_X\co M\to X$ via $r_{i,i}$, and
define the control map $ p^{i,i} \co   M^{i,i}\to  X^{i,i}$
to be the restriction of $ p^{i,i}_*$ to $ M^{i,i}$.
Note that the stereographic projection from $(v_i,-H)\in \br^N\times\br$ defines
a homeomorphism $X \to X^{i,i}$ sending $V_k$ and $L_k$ to
$V^{i,i}_k$ and $L^{i,i}_k$ respectively, since $V_i$ is the
cone on $L_i$ with center $v_i$.

Next, for each $i=1,\dots,m-1$,
define  $X^{i,i+1}_*\subset\br^N\times\br\times\br$ and its subsets $X^{i,i+1}$,
$V^{i,i+1}_k$, $L^{i,i+1}_k$  ($k=i+2,\cdots,m$) by
\begin{align*}
 X^{i,i+1}_*&= X^{i,i} \cup ( V^{i,i}_{i+1}\times [0,H])~,\\
 X^{i,i+1} &=( X^{i,i}- V^{i,i}_{i+1})\cup
( L^{i,i}_{i+1}\times [0,H]) \cup ( V^{i,i}_{i+1}\times \{H\})
\subset X^{i,i}\times\br~,\\
 V^{i,i+1}_k &=(V^{i,i}_k-V^{i,i}_i)\cup (V^{i,i}_k\cap L^{i,i}_i\times [0,H])
\cup (V^{i,i}_k\cap V^{i,i}_i\times \{H\})~,\\
 L^{i,i+1}_k &=(L^{i,i}_k-V^{i,i}_i)\cup (L^{i,i}_k\cap L^{i,i}_i\times [0,H])
\cup (L^{i,i}_k\cap V^{i,i}_i\times \{H\})~.
\end{align*}
Again we use the product metric of $\br^N\times\br$ and $\br$.
The projection $\br^N\times\br\times\br\to \br^N\times\br$ restricts to a retraction
$r_{i,i+1}\co X^{i,i+1}_*\to X^{i,i}$.
The control maps
$ p^{i,i+1}_* \co   M^{i,i+1}_*\to X^{i,i+1}_*$ and
$ p^{i,i+1}\co   M^{i,i+1}\to X^{i,i+1}$ are defined to be the pullbacks of
$p^{i,i}_*$ via $r_{i,i+1}$ and $r_{i,i+1}|X^{i,i+1}$, respectively.
Although $V^{i,i}_{i+1}$ is not a cone, it is homeomorphic to $V_{i+1}$
and has a topological cone structure.  So one can construct a homeomorphism
from $X^{i,i+1}$ to $X^{i,i}$ sending $V^{i,i+1}_k$ and $L^{i,i+1}_k$ to
$V^{i,i}$ and $L^{i,i}$ respectively, and hence a homeomorphism to $X$.

We can continue this to inductively obtain the metric space
\[
 X^{i,j}_*=X^{i,j-1}\cup ( V^{i,j-1}_{j}\times[0,H])
\]
as a subset of $\br^N\times \br^{j-i+1}$,
and its subsets $ X^{i,j}\supset  V^{i,j}_k \supset L^{i,j}_k$ ($k=j+1,\cdots,m$),
together with control maps $p^{i,j}_*\co  M^{i,j}_*\to  X^{i,j}_*$,
and $ p^{i,j}\co  M^{i,j}\to  X^{i,j}$.
Topologically all the spaces $ X^{i,j}$'s are equal to $X$, and
all the sets $ V^{i,j}_k$'s are equal to $V_k$.
We are only changing the metric, the cell structure of $X$, and the control map.

\begin{center}
\includegraphics{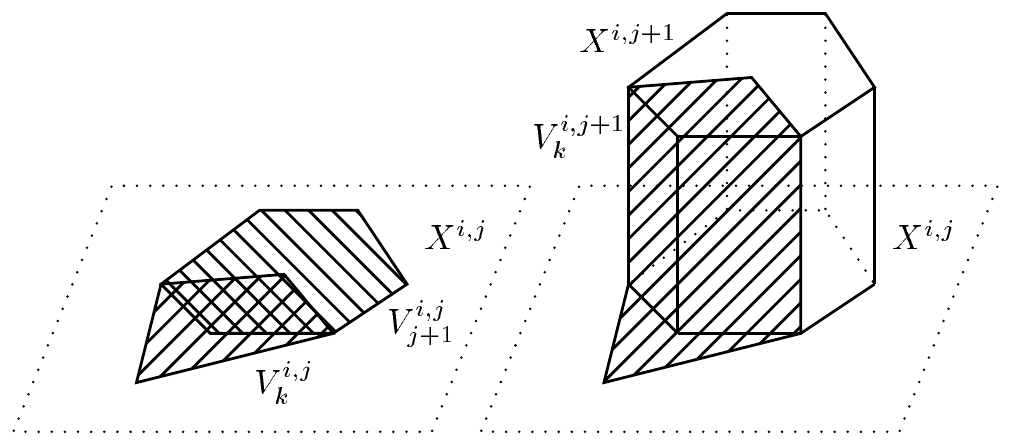}
\end{center}

Our next task is to do Alexander tricks at $v_1$, \dots, $v_m$ on these spaces
instead of $X$.
Since $\epsilon'\le\delta_0\le\delta$, we can split the original complex $c$ into
two pieces on $ V_1$ and the closure of its complements
by a $\lambda\epsilon'$ cobordism.
Now we construct a tower:
we make copies of the trivial cobordism from the pair on
$(V_1,  L_1)$ to itself and successively attach them to the cobordism along
$ V_1\times [0,H]$.
This is actually done on $ M^{1,1}_*$.

We use enough layers so that the result is a projective cobordism of
radius $\mu\lambda\epsilon'$ measured on $ X^{1,1}_*$
from $c=\wbar c_0$ to a complex $c'_1$ on $ p^{1,1}$.
Recall $\mu=100$ and it comes from taking a union of Poincar\'e cobordisms.
As described in previous sections, we can replace this by a free
cobordism of radius $\mu\lambda^2\epsilon'$ from $c$ to a free complex
$\wbar c_1$ on $ p^{1,1}$.

We postpone the squeezing to a later stage and go ahead to perform
Alexander trick over $V^{1,1}_2\subset X^{1,1}$ on $\wbar c_1$.
Although $X^{1,1}$ has a different metric from $X$,
the difference lies along the cylinder $L_1\times[0,H]$.
If $H$ is sufficiently large, then a rectification for $X^{1,1}$
can be easily constructed from those for $X$ and $[0,H]$,
and the $\delta$ and $\lambda$ for $X$ works also for $X^{1,1}$.
Since $\mu\lambda^2\epsilon'\le\delta$, we can do splitting and cut
out the portion on $V^{1,1}_2$ by a $\mu\lambda^3\epsilon'$ cobordism.
Again use enough copies of this to get a $\mu^2\lambda^3\epsilon'$ cobordism
on $p_*^{1,2}$ to a complex $\wbar c'_2$ on $p^{1,2}$ and then
replace this by free $\mu^2\lambda^4\epsilon'$ cobordism
to a free complex $\wbar c_2$ on $p^{1,2}$.
Since $\epsilon'\le\delta_0$, we can continue this process
to obtain a consecutive sequence of free cobordisms:
\begin{multline*}
\bfig
  \morphism/-/<800,0>[c=\wbar{c}_{0}{}_-X_*^{1,1}`
    \wbar{c}_1 {}_- X_*^{1,2};
    \mu\lambda^2\epsilon']
  \morphism(800,0)/-/<700,0>[\phantom{\wbar{c}_1 {}_- X_*^{1,2}}`
    \wbar{c}_2;
    (\mu\lambda^2)^2\epsilon']
  \morphism(1500,0)/-/<500,0>[\phantom{\wbar{c}_2}`\cdots;]
\efig \\
\bfig
  \morphism/-/<700,0>[\cdots`
    \wbar{c}_{m-2} {}_- X_*^{1,m-1};]
  \morphism(700,0)/-/<1100,0>[\phantom{\wbar{c}_{m-2} {}_- X_*^{1,m-1}}`
    \wbar{c}_{m-1} {}_- X_*^{1,m};
    (\mu\lambda^2)^{m-1}\epsilon']
  \morphism(1800,0)/-/<800,0>[\phantom{\wbar{c}_{m-1} {}_- X_*^{1,m}}`
    \wbar{c}_m;(\mu\lambda^2)^{m}\epsilon']
\efig
\end{multline*}
%\[
%\xymatrix@C+4pt{
%c=\bar c_0~ \ar@{.}[r]^-{\mu\lambda^2\epsilon'}_-{X_*^{1,1}} &
%~\bar c_1~ \ar@{.}[r]^-{(\mu\lambda^2)^2\epsilon'}_-{X_*^{1,2}} &
%~\bar c_2~ \ar@{.}[r] &
%\dots \ar@{.}[r] &
%~\bar c_{m-2}~\ar@{.}[r]^-{(\mu\lambda^2)^{m-1}\epsilon'}_-{X_*^{1,m-1}} &
%~\bar c_{m-1}~\ar@{.}[r]^-{(\mu\lambda^2)^{m}\epsilon'}_-{X_*^{1,m}} &
%~\bar c_m
%}
%\]
Now we construct a map $S^{1,m}\co X^{1,m}\to X$ and a map
$\wtilde{S}^{1,m}\co M^{1,m}\to M$ which covers $S^{1,m}$
so that the functorial image of $\wbar c_m$ has the desired property.
This is done by inductively constructing maps
$S^{i,j}_*\co X^{i,j}_*\to X$ and its restriction $S^{i,j}\co X^{i,j}\to X$
covered by maps $\wtilde{S}^{i,j}_*\co M^{i,j}_*\to M$ and
$\wtilde{S}^{i,j}\co M^{i,j}\to M$, respectively,  for certain pairs $j\ge i$.

First we define $S^{i,i}_*\co X^{i,i}_*\to X$.
Let us recall that $S_i\subset V_i$ denotes the stable set at $v_i$.
Using the strong deformation retraction $s_{i,t}\co V_i\to V_i$, define a map
$S'_i\co X^{i,i}_*\to X^{i,i}_*$ by
\[
(x,h) \mapsto \begin{cases} (x,0) &\hbox{if $x\in X$ and $h=0$,}\\
(s_{i,h/H}(x), h) & \hbox{if $x\in V_i$ and $h>0$}
\end{cases}
\]
This map is covered by a map $\wtilde{S}'_i\co M^{i,i}_*\to M^{i,i}_*$~.

\begin{lemma}
$S'_i$ has Lipschitz constant $K+1$.
\end{lemma}

\begin{proof}
This is obtained by applying \fullref{lipschitz} to the sets of the form
$\{x\}\cup V_i$ for $x\in X-V_i$, extending the map $s_i$ on $x$ by $s_i(x)=x$.
\end{proof}

$S^{i,i}_*\co X^{i,i}_*\to X$ is defined by composing $S'_i$ with the
projection $r_{i,i}\co X^{i,i}_*\to X$. It has Lipschitz constant $K+1$.
Since $r_{i,i}$ is obviously covered by a map $M^{i,i}_*\to M$,
the map $S^{i,i}_*$ is covered by a map $\wtilde{S}^{i,i}_*\co M^{i,i}_*\to M$.
Define $S^{i,i}\co X^{i,i}\to X$ to be the restriction of
$S^{i,i}_*$.

Now recall that $X^{1,2}_*$ and $X^{2,2}_*$ are obtained by attaching
$V^{1,1}_2\times [0,H]$ and $V_2\times [0,H]$ to $X^{1,1}$ and $X$, respectively.
Since $S^{1,1}\co X^{1,1}\to X$ maps $V^{1,1}_2$ to $V_2$, the product map
$S^{1,1}\times 1_{[0,H]}\co X^{1,1}\times [0,H]\to X\times [0,H]$ restricts to a map
$S^{1,1}\times 1|\co X^{1,2}_* \to X^{2,2}_*$.
Compose this with $S^{2,2}_*\co X^{2,2}_*\to X$ to define $S^{1,2}_*\co X^{1,2}\to X$
which is covered by a map $\wtilde{S}^{1,2}_*\co M^{1,2}_*\to M$.
Continue this process as in the following diagram to eventually get the desired map
$S^{1,m}\co X^{1,m}\to X$.
$$\bfig
  \barrsquare(0,400)|alra|/->`^{ (}->`^{ (}->`->/<800,400>[
    X^{1,j-1}`X`X^{1,j}_*`X^{j,j}_*;
    S^{1,j-1}```S^{1,j-1}\times1|]
  \barrsquare(0,0)|blrb|/{@{>}@/^-10pt/}`<-^)`=`->/<1600,400>[
    \phantom{X^{1,j}_*}`X`X^{1,j}`X;
    S^{1,j}_*```S^{1,j}]
  \morphism(800,400)<800,0>[\phantom{`X^{j,j}_*}`\phantom{X};S^{j,j}_*]
  \efig$$
Recall that there is a topological identification of $X^{1,m}$ with $X$.
So we can think of $S^{1,m}$ to be a map from $X$ to $X$ equipped with different
metrics.  Although it is not a homeomorphism, it preserves all the simplices,
{\it i.e.} $S^{1,m}(\Delta)=\Delta$ for every simplex $\Delta$ of $X$.
When restricted to a simplex, $S^{1,m}$ has Lipschitz constant
$(K+1)^m d/H$.

\begin{center}
\includegraphics{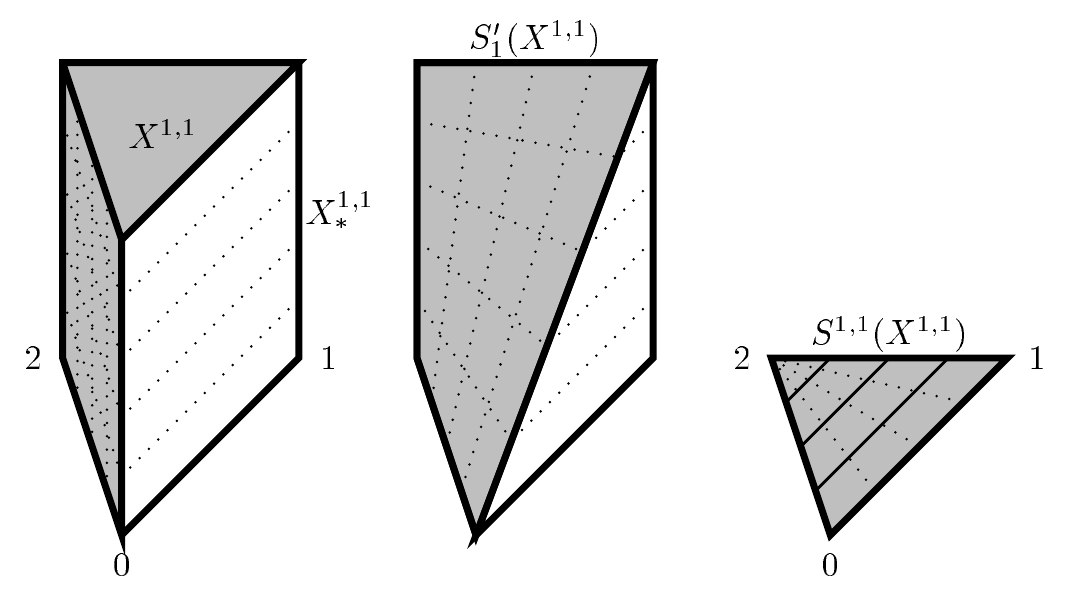}
\end{center}

The three pictures above illustrate the application of $S^{1,1}$ to $X^{1,1}$.
The thin solid lines in the rightmost picture indicate the direction
in which controls are obtained.

The three pictures below illustrate the application of $S^{1,2}$ to $X^{1,2}$.
The leftmost picture shows the image $(S^{1,1}\times 1)(X^{1,2})=X^{2,2}$.
Again the thin solid lines on the faces indicate the directions
in which controls are obtained.

\begin{center}
\includegraphics{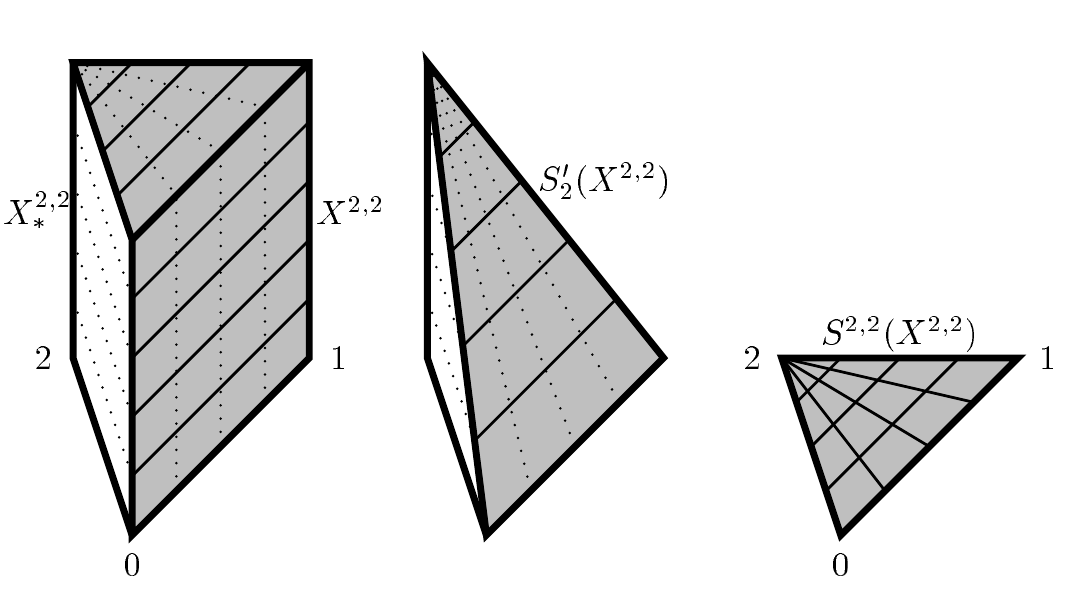}
\end{center}

Let us consider the functorial image $c_m$ of $\wbar c_m$ by the map
$\wtilde{S}^{1,m}\co M^{1,m}\to M$.
Recall that $\wbar c_m$ has radius $\epsilon''=(\mu\lambda^2)^m\epsilon'$.
Take a ball $B$ of radius $\epsilon''$ with in $X^{1,m}$.
$B$ is the union of subsets $B\cap\Delta$ each having diameter $2\epsilon''$,
where $\Delta$ are the simplices of $X^{1,m}$.
The images of $B\cap\Delta$ in $X$ by $S^{1,m}$ all have diameter
$2(K+1)^m\epsilon''d/H$ and their union $S^{1,m}(B)$ is connected.
Therefore $S^{1,m}(B)$ has diameter $2\sharp(X)(K+1)^m\epsilon''d/H$.
Thus $c_m$ has radius
\[
4\sharp(X)(K+1)^m(\mu\lambda^2)^m \epsilon'd/H~,
\]
and this is smaller than $\epsilon$ by the choice of $H$.

It remains to find a constant $\kappa$ such that $c$ and $c_m$ are
$\kappa\epsilon'$ cobordant.
Define a complex $c_i$ on $p_X$ to be the functorial image of $\wbar c_i$
by the map $\wtilde{S}^{1,i}\co M^{1,i}\to M$.
The functorial image of the $(\mu\lambda^2)^i\epsilon'$ cobordism between
$\wbar c_{i-1}$ and $\wbar c_i$ by the map $\wtilde{S}^{1,i}_*$
gives a $4\sharp(X)(K+1)^i(\mu\lambda^2)^i\epsilon'$ cobordism
between $c_{i-1}$ and $c_i$.
Composing these we get a $4\mu\sharp(X)(K+1)^m(\mu\lambda^2)^m\epsilon'$ cobordism
between $c$ and $c_m$.
Thus $\kappa=4\mu\sharp(X)(K+1)^m(\mu\lambda^2)^m$ works.
This completes the proof.
\end{proof}

Note that \fullref{squeeze} implies that the relax-control map in
\fullref{stability} is surjective:
Take an element $[c']\in L_n^{\delta',\epsilon'}(X;p_X,R)$ with
$\delta'\le\delta_0$.
Then the inequality $\epsilon'\le\delta_0$ holds and therefore
there is a Poincar\'e cobordism of radius $\kappa\epsilon'$ ($\le\delta_0$)
from $c'$ to a quadratic Poincar\'e complex $c$ of radius $\epsilon$,
determining an element  $[c]\in L_n^{\delta,\epsilon}(X;p_X,R)$
whose image under the relax-control map is $[c']$.

Squeezing for complexes can be generalized to squeezing for pairs.

\begin{lemma}[Squeezing of Quadratic Poincar\'e Pairs]
\label{squeeze2}
Let $n\ge 2$ be an integer and $X$ be a finite polyhedron.
There exist constants $\delta_0>0$ and $\kappa>1$
such that the following hold: If $\delta<\epsilon'\le\delta'\le\delta_0$, then
any $(n+1)$--dimensional quadratic Poincar\'e $R$--module pair of
radius $\delta'$ on a polyhedral stratified system of fibrations
over $X$ with $\epsilon'$ Poincar\'e boundary is $\kappa\delta'$ Poincar\'e
cobordant to a quadratic Poincar\'e pair of radius $\epsilon$.
The cobordism between the boundary is $\kappa\epsilon'$ Poincar\'e.
\end{lemma}

\begin{proof}
Same as the proof of \fullref{squeeze}. Use the Alexander trick
for pairs.
\end{proof}

\begin{corollary}[Relative Squeezing of Quadratic Poincar\'e Pairs]
\label{squeeze3}
Let $n\ge 2$ be an integer and $X$ be a finite polyhedron.
There exist constants $\delta_0>0$ and $\kappa>1$
such that the following hold: If $\kappa\epsilon<\delta'\le\delta_0$, then
any $(n+1)$--dimensional quadratic Poincar\'e $R$--module pair of
radius $\delta'$ on a polyhedral stratified system of fibrations
over $X$ with an $\epsilon$ Poincar\'e boundary
is $\kappa\delta'$ Poincar\'e cobordant fixing the boundary
to a quadratic Poincar\'e pair of radius $\kappa\epsilon$.
\end{corollary}

\begin{proof}
Temporarily choose $\delta_0$ and $\kappa$ as in \fullref{squeeze2}.
Suppose $\kappa\epsilon<\delta'\le\delta_0$, and
let $d=(f\co C\to D, (\delta\psi,\psi))$ be an $(n+1)$--dimensional
quadratic Poincar\'e pair of radius $\delta'$, and assume that
$(C,\psi)$ is $\epsilon$ Poincar\'e.
Choose a positive number $\epsilon'<\epsilon$.
By \fullref{squeeze2}, $d$ is $\kappa\delta'$ cobordant to
a quadratic Poincar\'e pair $d'=(f'\co C'\to D',(\delta\psi',\psi'))$
of radius $\epsilon$.
Glue $d'$ to the $\kappa\epsilon$ Poincar\'e cobordism between
$(C,\psi)$ and $(C',\psi')$ to get a quadratic Poincar\'e pair
$d''=(f''\co C\to D'', (\delta\psi'',\psi))$ of radius $100\kappa\epsilon$.
By construction, $d\cup -d''$ is $100\kappa\delta'$ Poincar\'e null-cobordant.
Thus $100\kappa$ works as the $\kappa$ in the statement of the lemma.
\end{proof}

The injectivity of the relax-control map follows from this:
Temporarily let $\delta_0$ and $\kappa$ be as in \fullref{squeeze3}, and
suppose $\delta\le\delta'$, $\epsilon\le\epsilon'$, $\kappa\epsilon\le\delta$.
Take an element $[c]$ in the kernel of the relax control map
\[
L_n^{\delta,\epsilon}(X;p_X,R)\to L_n^{\delta',\epsilon'}(X;p_X,R)~.
\]
By \fullref{nullcobordism}, 
the quadratic complex $c=(C,\psi)$ of radius $\epsilon'$ is the boundary of
an $(n+1)$--dimensional quadratic Poincar\'e pair $(f\co C\to D,(\delta\psi,\psi))$
of radius $100\delta'$.
If $\delta'\le \delta_0/100$, then
$\kappa\epsilon\le 100\delta'\le\delta_0$, and by \fullref{squeeze3}
the element $[c]$ is 0 in $L_n^{\kappa\epsilon,\epsilon}(X;p_X,R)$, and
hence also in $L_n^{\delta,\epsilon}(X;p_X,R)$.
So, by replacing $\delta_0$ with $\delta_0/100$, we established the desired
injectivity.
This finishes the proof of \fullref{stability}.

%%%%%%%%%%%%%%%%%%%%%%%%%%%%%%%%%%%%%%%%%%%%%%%%%%%%%%%%%%%%%%%%%%%%%%
\section{Variations}
\label{section_variations}

\subsection{Projective $L$--groups}
There is a controlled analogue of projective $L^p$--groups.
$L_n^{p,\delta,\epsilon}(X;p_X,R)$ is defined using
$\epsilon$ Poincar\'e $\epsilon$ quadratic projective $R$--module complexes on $p_X$
and $\delta$ Poincar\'e $\delta$ projective cobordisms.
Similar stability results hold for these.

To get a squeezing result in the $L^p$--group case, we first take the tensor product
of the given projective quadratic Poincar\'e complex $c$ with the symmetric complex
$\sigma(S^1)$ of the circle $S^1$.
By replacing it with a finite cover if necessary,
we may assume that the radius of $\sigma(S^1)$ is sufficiently small.
If the radius of $c$ is also sufficiently small,
we can construct a cobordism to a squeezed complex.
Split the cobordism along $X\times \{\hbox{two points}\}\subset X\times S^1$
to get a projective cobordism from the
original complex to a squeezed projective complex.

\subsection{$UV^1$ control maps}
When the control map is $UV^1$, there is no need to use paths to define
morphisms between geometric modules (see Pedersen--Quinn--Ranicki
\cite{PQR}).  This simplifies the situation
quite a lot, and we have:

\begin{proposition}
Let $p_X\co M\to X$ be a $UV^1$ map to a finite polyhedron. Then for
any pair of positive numbers $\delta\ge\epsilon$, there is an
isomorphism
\[
L_n^{\delta,\epsilon}(X;p_X,R)\cong L_n^{\delta,\epsilon}(X;1_X,R)
\]
for any ring with involution $R$ and any integer $n\ge 0$.
\end{proposition}

By \fullref{stability}, the stability holds for
$L_n^{\delta,\epsilon}(X;1_X,R)$ and hence the stability holds also for
$L_n^{\delta,\epsilon}(X;p_X,R)$.

\subsection{Compact metric ANR's}
Squeezing and stability also hold when $X$ is a compact metric ANR, and
the control map is a fibration.
To see this, embed $X$ in the Hilbert cube $I^\infty$.
There is a closed neighborhood $U$ of $X$ of the form $P\times I^{\infty - N}$,
where $P$ is a polyhedron in $I^N$.
Use the fact that the retraction from $U$ to $X$ is uniformly continuous
to deduce the desired stability from the stability on $P$ and $U$.

%%%%%%%%%%%%%%%
\section{Relations to bounded $L$--theory}
\label{boundedL}
In this section we shall identify the controlled
$L$--theory groups with a bounded $L$--theory group, at least in the
case of constant coefficients. The main advantage to having a
bounded controlled description, is that it facilitates
computations.

\begin{definition} Let $X$ be a finite polyhedron and $R$ a ring with
involution. Let $p_X\co X\times K\to X$ be a trivial fibration. We denote
the common value of $L_n^{\delta,\epsilon}(X;p_X,R)$ for small values of
$\delta$ and $\epsilon$, which exists by \fullref{stability}, by $L_n^{h,c}(X;p_X,R)$.
Here the $h$ signifies that we have no simpleness condition and the $c$ stands
for controlled.
\end{definition}

We may embed the finite polyhedron $X$ in a large dimensional sphere $S^n$ and
consider the open cone $O(X)=\{t\cdot x \in R^{n+1}| t\in [0,\infty), x\in X\}$.
We denote $X$ with a disjoint basepoint added by $X_+$.

\begin{theorem} Let $p_X\co X\times K \to X$ be as above, $\pi=\pi_1(K)$, $R$ a
ring with involution. Then
\[
L_n^{c,h}(X;p_X,R)\cong L_{n+1}^s(\mathcal C_{O(X_+)}(R[\pi]))
\]
where $\mathcal C_{O(X_+)}(R[\pi])$ denotes the category of free $R[\pi]$
modules parameterized by $O(X_+)$ and bounded morphisms.
\end{theorem}

\begin{proof} Given an element in $L_n^{c,h}(X;p_X,R)$, we can choose a stable
$(\delta,\epsilon)$ representative. Crossing with the symmetric
chain complex of $(-\infty,0]$ produces a bounded quadratic chain
complex when we parameterize it by $O(+)$, which is obviously a
half line, with $+$ being the extra basepoint. According to
\fullref{stability}, we may produce a sequence of bordisms to
increasingly smaller representatives of the given element in
$L_n^{c,h}(X;p_X,R)$. These bordisms may be parameterized by
$\{t\cdot x|x\in X, a_i<t<a_{i+1}\}$ where the sequence of $a_i$'s
is chosen such that when these bordisms are glued together, we
obtain a bounded quadratic complex parameterized by $O(X_+)$. We
get an $s$--decoration because obviously we can split the bounded
quadratic complex. The map in the opposite direction is given by a
splitting obtained the same way as in \fullref{splitting_lemma}.
\end{proof}

One advantage of a categorical description is computational. We have as close
an analogue to excision as is possible in the following: Let $Y$ be a
subcomplex of $X$, and $S$ a ring with involution. We then get a sequence of
categories
\[
\mathcal C_{O(Y_+)}(S)\to \mathcal C_{O(X_+)}(S)\to \mathcal C_{O(X/Y)}(S)
\]
which leads to a long exact sequence
\[
\ldots \to L_n^a(\mathcal C_{O(Y_+)}(S))\to L_n^b(\mathcal C_{O(X_+)}(S))\to
L_n^c(\mathcal C_{O(X/Y)}(S))\to \ldots
\]
where the rule to determine the decorations is that $b$ can be
chosen to be any involution preserving subgroup of $K_i(\mathcal
C_{O(X_+)}(S))$, $i\le 2$, but then $c$ has to be the image in
$K_i(\mathcal C_{O(X/Y)}(S))$, and $a$ has to be the preimage in
$K_i(\mathcal C_{O(Y_+)}(S))$. See the paper \cite{hpII} by Hambleton
and Pedersen for a derivation of these exact sequences. This makes it
possible to do extensive calculations with controlled $L$--groups.

\bibliographystyle{gtart}
\bibliography{link}

\end{document}